\documentclass[10pt,twoside]{amsart}

\usepackage{amsthm}
\usepackage{amsmath}
\usepackage{amssymb}
\usepackage{amsfonts}
\usepackage{latexsym}
\usepackage{enumitem}
\usepackage{mathrsfs}
\usepackage[all]{xy} \SelectTips{eu}{}
\usepackage{hyperref}

\newtheorem*{intsetup*}{Setup}

\newtheorem*{intNotations*}{Notations and Conventions}

\newtheorem{intthm}{Theorem}[]

\newtheorem*{intque*}{Question}
\newtheorem*{intexa*}{Example}

\newcommand{\numberseries}{\bfseries}   
\newlength{\thmtopspace}                
\newlength{\thmbotspace}                
\newlength{\thmheadspace}               
\newlength{\thmindent}                  
\newtheoremstyle{bfupright head,slanted body}
                {\thmtopspace}{\thmbotspace}
                {\slshape}{\thmindent}{\bfseries}{.}{\thmheadspace}
                {{\numberseries \thmnumber{#2\;}}\thmnote{#3}}

\newtheoremstyle{bfupright head,upright body}
                {\thmtopspace}{\thmbotspace}
                {\upshape}{\thmindent}{\bfseries}{.}{\thmheadspace}
                {{\numberseries \thmnumber{#2\;}}\thmnote{#3}}

\newtheoremstyle{fixed bf head,slanted body}
                {\thmtopspace}{\thmbotspace}{\slshape}
                {\thmindent}{\bfseries}{.}{\thmheadspace}
                {{\numberseries \thmnumber{#2\;}}\thmname{#1}\thmnote{ (#3)}}

\newtheoremstyle{fixed bf head,upright body}
                {\thmtopspace}{\thmbotspace}{\upshape}
                {\thmindent}{\bfseries}{.}{\thmheadspace}
                {{\numberseries \thmnumber{#2\;}}\thmname{#1}\thmnote{ (#3)}}

\newtheoremstyle{numbered paragraph}
                {\thmtopspace}{\thmbotspace}{\upshape}
                {\thmindent}{\upshape}{}{\thmheadspace}
                {{\numberseries \thmnumber{#2.}}}

\theoremstyle{bfupright head,slanted body}
\newtheorem{res}{}[section]             \newtheorem*{res*}{}

\theoremstyle{bfupright head,upright body}
\newtheorem{bfhpg}[res]{}               \newtheorem*{bfhpg*}{}

\theoremstyle{fixed bf head,slanted body}
          \newtheorem*{thm*}{Theorem}
\newtheorem{prp}[res]{Proposition}      \newtheorem*{prp*}{Proposition}
\newtheorem{cor}[res]{Corollary}        \newtheorem*{cor*}{Corollary}
\newtheorem{lem}[res]{Lemma}            \newtheorem*{lem*}{Lemma}
         \newtheorem*{que*}{Question}
\theoremstyle{fixed bf head,upright body}
\newtheorem{setup}[res]{Setup}           \newtheorem*{setup*}{Setup}
\newtheorem{dfn}[res]{Definition}       \newtheorem*{dfn*}{Definition}
\newtheorem{rmk}[res]{Remark}           \newtheorem*{rmk*}{Remark}
\newtheorem{exa}[res]{Example}           \newtheorem*{exa*}{Example}

\theoremstyle{numbered paragraph}
\newtheorem{ipg}[res]{}

\newlength{\thmlistleft}        
\newlength{\thmlistright}       
\newlength{\thmlistpartopsep}   
\newlength{\thmlisttopsep}      
\newlength{\thmlistparsep}      
\newlength{\thmlistitemsep}     

\newcounter{eqc}
  {\end{list}}%

\newcounter{prt}
\newenvironment{prt}{\begin{list}{\upshape (\alph{prt})}%
    {\usecounter{prt}%
      \setlength{\leftmargin}{\thmlistleft}%
      \setlength{\labelwidth}{\thmlistleft}%
      \setlength{\rightmargin}{\thmlistright}%
      \setlength{\partopsep}{\thmlistpartopsep}%
      \setlength{\topsep}{\thmlisttopsep}%
      \setlength{\parsep}{\thmlistparsep}%
      \setlength{\itemsep}{\thmlistitemsep}}}%
  {\end{list}}%

\newcounter{rqm}
  {\end{list}}%

  {\end{list}}%

\newenvironment{prf*}[1][Proof]{%
  \begin{proof}[\bf #1]
    \setcounter{equation}{0}
    }
  {\end{proof}
}

\newcommand{\pgref}[1]{\ref{#1}}

\renewcommand{\eqref}[1]{(\pgref{eq:#1})}

\newcommand{\thmcite}[2][?]{\cite[Theorem~#1]{#2}}
\newcommand{\prpcite}[2][?]{\cite[Proposition~#1]{#2}}

\newcommand{\corcite}[2][?]{\cite[Corollary~#1]{#2}}
\newcommand{\lemcite}[2][?]{\cite[Lemma~#1]{#2}}

\newcommand{\exacite}[2][?]{\cite[Example~#1]{#2}}

\numberwithin{equation}{res}

\def\urltilda{\kern -.15em\lower .7ex\hbox{\~{}}\kern .04em}

\newcommand{\GF}{\mathsf{GF}}

\newcommand{\PGF}{\mathsf{PGF}}

\newcommand{\sfT}{\mathsf{T}}

\newcommand{\Proj}{\mathsf{Proj}}
\newcommand{\id}{\mathrm{id}}
\newcommand{\fd}{\mathrm{fd}}
\newcommand{\pd}{\mathrm{pd}}
\newcommand{\Inj}{\mathsf{Inj}}
\newcommand{\Flat}{\mathsf{Flat}}
\newcommand{\Mod}{\mathsf{Mod}}

\newcommand{\Ind}{\mathrm{Ind}}
\newcommand{\Id}{\mathrm{Id}}
\newcommand{\Rop}{R^{\sf op}}
\newcommand{\op}{\sf op}

\newcommand{\cok}{\mbox{\rm coker}}
\newcommand{\kernel}{\mbox{\rm ker}}
\newcommand{\im}{\mbox{\rm Im}}

\newcommand{\Coind}{\mathrm{Coind}}

\newcommand{\Hom}{\operatorname{Hom}}

\newcommand{\Tor}{\operatorname{Tor}}

\newcommand{\is}{\cong}

   \def\soft#1{\leavevmode\setbox0=\hbox{h}\dimen7=\ht0\advance
    \dimen7 by-1ex\relax\if t#1\relax\rlap{\raise.6\dimen7
    \hbox{\kern.3ex\char'47}}#1\relax\else\if T#1\relax
    \rlap{\raise.5\dimen7\hbox{\kern1.3ex\char'47}}#1\relax
    \else\if d#1\relax\rlap{\raise.5\dimen7\hbox{\kern.9ex
    \char'47}}#1\relax\else\if D#1\relax\rlap{\raise.5\dimen7
    \hbox{\kern1.4ex\char'47}}#1\relax\else\if l#1\relax
    \rlap{\raise.5\dimen7\hbox{\kern.4ex\char'47}}#1\relax
    \else\if L#1\relax\rlap{\raise.5\dimen7\hbox{\kern.7ex
    \char'47}}#1\relax\else\message{accent \string\soft
    \space #1 not defined!}#1\relax\fi\fi\fi\fi\fi\fi}

\begin{document}

\title[(projectively coresolved) gorenstein flat modules over tensor rings]%
{(projectively coresolved) gorenstein flat modules over tensor rings}

\author[G.L. Tang]{Guoliang Tang}
\address{Guoliang Tang: School of Mathematical Sciences,
Zhejiang Normal University, Jinhua 321004, China}
\email{tangguoliang970125@163.com}

\author[J.Q. Wei]{Jiaqun Wei}
\address{Jiaqun Wei (Corresponding author): School of Mathematical Sciences,
Zhejiang Normal University, Jinhua 321004, China}
\email{weijiaqun5479@zjnu.edu.cn}

\thanks{J.Q. Wei was partly supported by the National Natural Science
Foundation of China (Grant Nos. 12571042, 12271249)
and the Natural Science Foundation of Zhejiang Province (Grant No. LZ25A010002)}


\keywords{Gorenstein flat module;
projectively coresolved Gorenstein flat module; tensor ring;
trivial ring extension; Morita context ring.}

\footnotetext{2020 \emph{Mathematics Subject Classification}. 18G25; 16D90.}


\begin{abstract}
Let $T_R(M)$ be a tensor ring, where $R$ is a ring and
$M$ is an $N$-nilpotent $R$-bimodule.
Under certain conditions, we characterize projectively coresolved Gorenstein flat
modules over $T_R(M)$, showing that a $T_R(M)$-module $(X,u)$ is
projectively coresolved Gorenstein flat if and only if $u$ is monomorphic and
$\cok(u)$ is a projectively coresolved Gorenstein flat $R$-module.
A class of Gorenstein flat modules over $T_R(M)$ are also explicitly described.
We discuss applications to trivial ring extensions and Morita context rings.
\end{abstract}

\maketitle

\thispagestyle{empty}

\section*{Introduction}
\label{Preliminaries}
\noindent
Throughout the paper, all rings are assumed to be nonzero associative rings with identity, and all modules are assumed to be unitary. For a ring $R$, we adopt the convention that an $R$-module refers to a left $R$-module, while right $R$-modules are viewed as modules over the opposite ring $\Rop$.

Gorenstein homological algebra provides a powerful extension of classical homological techniques to broader algebraic contexts.
Within this framework, the notion of projectively coresolved Gorenstein flat modules (PGF, for short), introduced by \v{S}aroch and \v{S}\v{t}ov\'{\i}\v{c}ek \cite{GFEC}, refines the classical concepts of
Gorenstein projective modules of Enochs and Jenda \cite{1995EnochsGP} and
Gorenstein flat modules of Enochs, Jenda and Torrecillas \cite{GF1993}.
PGF modules are characterized by admitting projective coresolutions that remain exact after tensoring with every injective module,
and they have become an essential tool in constructing complete cotorsion pairs and in understanding definability and singular compactness phenomena.
Their foundational significance in modern development of relative homological algebra is now well established; see, for instance, \cite{2024PGFdim, 2020Alina, 2025PGF}.

Tensor rings provide a natural framework for studying homological invariants in the context of algebraic extensions.
Following Cohn's work on formal tensor algebras \cite{1991Cohn}, a tensor ring of an $R$-bimodule $M$ is defined by
\[
T_R(M)=\bigoplus_{i=0}^\infty M^{\otimes_R i},
\]
where $M^{\otimes_R0}=R$ and $M^{\otimes_R(i+1)}=M\otimes_R(M^{\otimes_Ri})$. This construction encodes iterated bimodule actions and includes as special cases trivial extensions and Morita context rings. When $M$ is $N$-\emph{nilpotent} with $N\geqslant 0$, that is, $M^{\otimes_R(N+1)}=0$, the tensor ring acquires a useful structural feature: as shown by Chen and Lu \cite{TRchen}, every $T_R(M)$-module can be uniquely represented by a pair $(X,u)$, where $X$ is an $R$-module and $u\colon M\otimes_R X\to X$ is an $R$-homomorphism.

Working over a Noetherian ring $R$ and an $N$-nilpotent $R$-bimodule $M$
that is finitely generated on both sides,
Chen and Lu \cite{TRchen} characterized the finitely generated Gorenstein projective
modules over the tensor ring $T_R(M)$.
Their theorem states that, under certain conditions, a
$T_R(M)$-module $(X,u)$ is finitely generated Gorenstein projective if and only if
$u$ is a monomorphism and its cokernel is a finitely generated
Gorenstein projective $R$-module.
This established tensor rings as a flexible and powerful framework in Gorenstein homological algebra.
The recent work of Di et al. \cite{2025DLST} generalizes this result to an arbitrary associative ring $R$ and an $R$-bimodule $M$ that is not necessarily finitely generated. This development naturally motivates the central question of this work.

\begin{que*}
How can one describe (projectively coresolved) Gorenstein flat modules over $T_R(M)$ in terms of $R$-modules, when $M$ is an $N$-nilpotent $R$-bimodule?
\end{que*}

The main purpose of the paper is to answer this question and to provide further applications. Throughout the rest of this section, we \textbf{assume} that $M$ is an $N$-nilpotent $R$-bimodule.
For a subcategory $\sf{X}$ of $R$-modules, we write $\Phi(\sf{X})$ for the subcategory of the category of $T_R(M)$-modules defined by
\begin{align*}
\Phi({\sf{X}}) = \{(X,u) \in \Mod(T_R(M)) \,|\,
u \textrm{ is a monomorphism and } \cok(u) \in \sf{X}\}.
\end{align*}
For the $R$-bimodule $M$, we consider the following Tor-vanishing condition:
\[
\tag{$\sfT$}
\Tor_{\geqslant 1}^R(M,\, M^{\otimes_R i}\otimes_R P)=0
\quad\text{for each } P\in\Proj(R) \text{ and } i\geqslant 1,
\]
which is inspired by the work of Chen and Lu \cite{TRchen} and plays a crucial role in establishing our main results. Under this condition, our results can be stated as follows; unexplained notation is collected in Section~\ref{Preliminaries}.

\begin{intthm}\label{THM PGF com}
Let the $R$-bimodule $M$ satisfy condition $(\sfT)$.
If $\pd_R M<\infty$ and $\fd_{\Rop} M<\infty$, then
$\PGF(T_R(M))=\Phi(\PGF(R))$.
\end{intthm}

Previous results, notably by Di et al. \cite{2025DLST},
have described Gorenstein flat $T_R(M)$-modules under coherence conditions.
We seek to extend this understanding to non-coherent rings.
Given that the characterization of flat $T_R(M)$-modules is inherently less refined
than that of projective modules (see Lemma \ref{projectives in F-Rep}),
our investigation naturally focuses on a salient and tractable subclass,
for which a complete characterization is provided.

\begin{intthm}\label{THM GF com}
Let the $R$-bimodule $M$ satisfy condition $(\sfT)$.
If $\fd_R M<\infty$ and $\fd_{\Rop} M<\infty$, then an $R$-module $X$ is Gorenstein flat if and only if $\Ind(X)$ is a Gorenstein flat $T_R(M)$-module.
\end{intthm}

We conclude the paper by presenting two applications of these theorems.
First, we investigate (projectively coresolved) Gorenstein flat modules over trivial extensions of rings, and subsequently we determine when a module over a Morita context ring is (projectively coresolved) Gorenstein flat.

The paper is organized as follows. Section~1 collects basic material on tensor rings and introduces notation used throughout the paper. Section~2 is devoted to the proof of Theorem~\ref{THM PGF com}. Section~3 contains the proof of Theorem~\ref{THM GF com}. Finally, Section~4 presents the applications.
\section{Preliminaries}
\label{Preliminaries}
\noindent
In this section we fix our notation, recall relevant concepts, and present several basic facts used throughout the paper. For a ring $R$, we denote by $\Mod(R)$ the category of left $R$-modules, and by $\Proj(R)$, $\Inj(R)$, and $\Flat(R)$ the full subcategories of $\Mod(R)$ consisting of projective, injective, and flat $R$-modules, respectively. For an $R$-module $X$, we write $\pd_R X$, $\id_R X$, and $\fd_R X$ for its projective, injective, and flat dimensions, and we denote by $X^{+}$ the character module $\Hom_{\mathbb{Z}}(X,\mathbb{Q}/\mathbb{Z})$.

\begin{bfhpg}[\bf Tensor rings]\label{Tensor rings}
Let $R$ be a ring and $M$ an $R$-bimodule. We set $M^{\otimes_R 0}=R$ and
$M^{\otimes_R(i+1)}=M\otimes_R(M^{\otimes_R i})$ for $i\geqslant 0$.
For an integer $N\geqslant 0$,
the $R$-bimodule $M$ is said to be $N$-\emph{nilpotent} if $M^{\otimes_R(N+1)}=0$.
Whenever $M$ is $N$-nilpotent, the \emph{tensor ring} with respect to $M$ is defined by
\[
T_R(M)=\bigoplus_{i=0}^N M^{\otimes_R i}.
\]

By \cite{TRchen}, the category $\Mod(T_R(M))$ is equivalent to the category $\Gamma$ whose objects are pairs $(X,u)$, where $X\in\Mod(R)$ and
$u\colon M\otimes_R X\to X$ is an $R$-homomorphism; a morphism
$f:(X,u)\to (X',u')$ in $\Gamma$ is an $R$-homomorphism $f\colon X\to X'$ satisfying
$f\circ u = u'\circ (M\otimes f)$, that is
$$\xymatrix{
  M\otimes_R X \ar[d]_{u} \ar[r]^{M\otimes f}
                & M\otimes_R X' \ar[d]^{u'}  \\
  X  \ar[r]^{f}
                & \, X'.             }$$
Similarly, the category $\Mod(T_R(M)^{\op})$ is equivalent to the category $\Omega$ consisting of pairs $[Y,v]$ with
$Y\in\Mod(\Rop)$ and $v\colon Y\to\Hom_{\Rop}(M,Y)$ an $\Rop$-homomorphism, with morphisms defined analogously.

Unless stated otherwise, a $T_R(M)$-module is always viewed as a pair $(X,u)$ with
$X\in\Mod(R)$ and $u\in\Hom_R(M\otimes_R X,X)$, and a $T_R(M)^{\op}$-module is
considered as a pair $[Y,v]$ with
$Y\in\Mod(\Rop)$ and $v\in\Hom_{\Rop}(Y,\Hom_{\Rop}(M,Y))$.
A sequence
\[
(X,u)\xrightarrow{f}(X',u')\xrightarrow{g}(X'',u'')
\]
in $\Mod(T_R(M))$ is exact if and only if the underlying sequence
$X\xrightarrow{f}X'\xrightarrow{g}X''$ is exact in $\Mod(R)$.
The analogous statement holds for modules over $T_R(M)^{\op}$.
\end{bfhpg}

\setup
Throughout the remainder of the paper, $M$ always denotes an $N$-nilpotent $R$-bimodule.

\begin{ipg}\label{adjoint pairs}
We recall from \cite{TRchen} the following adjoint pairs of functors:
\[
\xymatrix@C=3pc{
\Mod(R)
\ar[r]^-{S}
&
\Mod(T_R(M))
\ar@/_1.8pc/[l]_-{C}
\ar[r]^-{U}
&
\Mod(R)
\ar@/_1.8pc/[l]_-{\Ind}.
}
\]
The \emph{stalk functor} $S$ sends an $R$-module $X$ to $(X,0)$. Its left adjoint $C$ is given by
$\!C((X,u))=\cok(u)$; for a morphism
$f:(X,u)\to(Y,v)$, the map $C(f):\cok(u)\to\cok(v)$ is induced from $f$.
The \emph{forgetful functor} $U$ sends $(X,u)$ to $X$, and its left adjoint $\Ind$ is defined by
\[
\Ind(X)=\left(\bigoplus_{i=0}^N M^{\otimes_R i}\otimes_R X,\ c_X\right),
\]
where $c_X$ is the canonical inclusion
\[
M\otimes_R\left(\bigoplus_{i=0}^N M^{\otimes_R i}\otimes_R X\right)
\cong \bigoplus_{i=1}^N M^{\otimes_R i}\otimes_R X
\longrightarrow \bigoplus_{i=0}^N M^{\otimes_R i}\otimes_R X.
\]
For an $R$-homomorphism $f\colon X\to Y$, the morphism
$\Ind(f)\colon\Ind(X)\to\Ind(Y)$ is represented by the diagonal matrix whose entries are
$M^{\otimes_R i}\otimes_R f$ for $0\le i\le N$.
It is clear that $C\circ\Ind=\Id_{\Mod(R)}$, and by the Eilenberg--Watts theorem $\Ind$ is naturally isomorphic to the tensor functor $T_R(M)\otimes_R -$.
\end{ipg}

\begin{ipg}\label{adjoint pairs-dual}
Dually to \ref{adjoint pairs}, there are adjoint pairs
\[
\xymatrix@C=3pc{
\Mod(\Rop)
\ar[r]^-{S}
&
\Mod(T_R(M)^{\op})
\ar@/^1.8pc/[l]_-{K}
\ar[r]^-{U}
&
\Mod(\Rop)
\ar@/^1.8pc/[l]_-{\Coind}.
}
\]
Here $K([Y,v])=\ker(v)$, and for $Y\in\Mod(\Rop)$ we set
\[
\Coind(Y)=\left[\bigoplus_{i=0}^N\Hom_{\Rop}(M^{\otimes_R i},Y),\ r_Y\right],
\]
where $r_Y$ is the canonical map whose matrix form is
\[
\begin{pmatrix}
0 & 1 & 0 & \cdots & 0 \\
0 & 0 & 1 & \cdots & 0 \\
\vdots & \vdots & \vdots & \ddots & \vdots \\
0 & 0 & 0 & \cdots & 1
\end{pmatrix}.
\]
\end{ipg}

We now record several auxiliary results that will be used frequently.

\begin{lem} \label{projectives in F-Rep} \lemcite[1.9]{2025DLST}
Let $(X,u)$ be a $T_R(M)$-module and $[Y,v]$ a $T_R(M)^{\op}$-module.
Then the following statements hold.
\begin{prt}
\item
$(X,u)$ is projective if and only if there exists a projective $R$-module $P$
such that $(X,u) \cong \Ind(P)$.

\item
$[Y,v]$ is injective if and only if there exists an injective $\Rop$-module $E$
such that $[Y,v] \cong \Coind(E)$.

\item
$(X,u)$ is flat if and only if
$u$ is a monomorphism and $\cok(u)$ is a flat $R$-module.
\end{prt}
\end{lem}

\begin{lem} \label{lemma4.5 for (T)} \lemcite[2.6]{2025DLST}
Let $X$ be an $R$-bimodule satisfying the condition $(\sfT)$.
If $\fd_RX$ (resp., $\fd_{\Rop}X$) is finite,
then so is $\fd_RX^{\otimes_Ri}$ (resp., $\fd_{\Rop}X^{\otimes_Ri}$)
for each $i \geqslant 1$.
\end{lem}

The next result describes tensor products of $T_R(M)$-modules.

\begin{lem} \label{tensor products of TRM modules}
Let $(X,u)$ be a $T_R(M)$-module and $[Y,v]$ a $T_R(M)^{\op}$-module.
Then there is an isomorphism
$$[Y,v]\otimes_{T_R(M)}(X,u)\cong (Y\otimes_RX)/H,$$
of abelian groups,
where the subgroup $H$ is generated by all elements of the form
\begin{align*}
v(y)(m)\otimes x\,\,\textrm{and}\,\, y \otimes u(m\otimes x)
\end{align*}
with $x \in X$, $y \in Y$ and $m \in M$.
\end{lem}
\begin{prf*}
Let $F$ be the free abelian group with basis $Y \times X$; that is,
$$F=\{\sum_{(y,x) \in Y \times X} n_{(y,x)}(y,x)\,\,|\,\,
y \in Y, x \in X \,\,\text{and}\,\, n_{(y,x)} \in \mathbb{Z} \}.$$
Then, by the definition of tensor products, $Y\otimes_RX=F/H_1$,
where $H_1$ is the subgroup of $F$ generated by all elements of
the following three forms:
\begin{align*}
(y_1+y_2,x)-(y_1,x)-(y_2,x); \\
(y,x_1+x_2)-(y,x_1)-(y,x_2); \\
(yr,x)-(y,rx) \qquad \quad \,\,
\end{align*}
with $y$, $y_1$, $y_2 \in Y$, $x$, $x_1$, $x_2 \in X$ and $r \in R$.
Similarly, by the definition of tensor products, we have
$[Y,v]\otimes_{T_R(M)}(X,u)=F/H_2$,
where $H_2$ is the subgroup of $F$ generated by all elements of
the following three types:
\begin{align*}
(y_1+y_2,x)-(y_1,x)-(y_2,x); \qquad \qquad \qquad \\
(y,x_1+x_2)-(y,x_1)-(y,x_2); \qquad \qquad \qquad \\
(y(r,m_{11},m_{21}\otimes m_{22},\cdots, m_{N1}\otimes m_{N2}\otimes \cdots \otimes m_{NN}),x)   \\
-(y,(r,m_{11},m_{21}\otimes m_{22},\cdots, m_{N1}\otimes m_{N2}\otimes \cdots \otimes m_{NN})x)
\end{align*}
with $y$, $y_1$, $y_2 \in Y$, $x$, $x_1$, $x_2 \in X$, $r \in R$
and $m_{ij} \in M$ for $j \leqslant i \leqslant N$.
Clearly, $H_1 \subseteq H_2$.
We have the isomorphism
$$[Y,v]\otimes_{T_R(M)}(X,u)=F/H_2 \cong
(F/H_1)/(H_2/H_1)=(X\otimes_RY)/(H_2/H_1).$$
Note that the factor group $H_2/H_1$ is generated by the images
of the generators of $H_2$, i.e., elements of the form
$$((v(y)(m), x)-(y, u(m\otimes x)))+H_1
=: v(y)(m)\otimes x-y\otimes u(m\otimes x)$$
with $y \in Y$, $x \in X$ and $m \in M$.
This completes the proof.
\end{prf*}

The next result is an immediate consequence of
Lemma \ref{tensor products of TRM modules}.

\begin{cor} \label{cor fof tensor products}
Let $(X,u)$ be a $T_R(M)$-module and $W$ a $\Rop$-module.
Then we have
$S(W)\otimes_{T_R(M)}(X,u)\cong (W\otimes_RX)/(W\otimes_R\im(u))
\cong W\otimes_R\cok(u)$
\end{cor}

\section{Proof of Theorem \ref{THM PGF com}}
\noindent
This section is dedicated to proving~\ref{THM PGF com}.
Throughout, we assume that the $R$-bimodule $M$ satisfies condition $(\sfT)$.
We begin by recalling the definition of projectively coresolved Gorenstein flat modules.

\begin{dfn}\label{def of silting} \cite{GFEC}
An $R$-module $X$ is said to be \emph{projectively coresolved Gorenstein flat}
(PGF, for short) if there exists an exact complex
\[
\xymatrix@C=0.5cm{
P^\bullet:\ \cdots \to P^{-1} \ar[r]^{\qquad \quad d^{-1}} & P^0
\ar[r]^{d^0 \quad \,\,} & P^1 \to \cdots
}
\]
of projective $R$-modules with $X \cong \ker(P^0 \to P^1)$ such that the complex
$E\otimes_R P^\bullet$ is exact for every injective $\Rop$-module $E$.
We denote by $\PGF(R)$ the class of projectively coresolved Gorenstein flat
modules over $R$.
\end{dfn}

\begin{rmk}\label{id finite}
According to \lemcite[2.1]{2025PGF}, an $R$-module $X$ is PGF if and only if there exists an exact complex $P^\bullet$ of projectives as above such that $E\otimes_R P^\bullet$ remains exact for every $\Rop$-module $E$ with finite injective dimension.
\end{rmk}

\begin{lem}\label{dim of S(R)}
The following assertions hold.
\begin{prt}
\item If $\fd_R M<\infty$, then $\fd_{T_R(M)} S(R)<\infty$.
\item If $\fd_{\Rop} M<\infty$, then $\fd_{T_R(M)^{\op}} S(R)<\infty$.
\end{prt}
\end{lem}
\begin{prf*}
We will prove (a); the proof of (b) is analogous.
Choose a flat resolution
\[
0 \to F^n \to \cdots \to F^1 \to F^0 \to M \to 0
\]
of $M$ as a left $R$-module. Since $\Tor_{\geqslant 1}^R(M, M^{\otimes_R i})=0$ for all $i\geqslant 1$ (by $(\sfT)$) and $M^{\otimes_R(N+1)}=0$,
we obtain the following exact sequences
\begin{align*}
0 \to M\otimes_RF^n \to \cdots \to M\otimes_RF^0 \to M\otimes_RM \to 0,
\qquad \qquad \,\,\,\, \\
0 \to M\otimes_RM\otimes_RF^n \to \cdots
\to M\otimes_RM\otimes_RF^0 \to M\otimes_RM\otimes_RM \to 0, \\
\cdots  \qquad \qquad \qquad \qquad \qquad \qquad \qquad \qquad \\
0 \to M^{\otimes_R(N-1)}\otimes_RF^n \to \cdots \to
M^{\otimes_R(N-1)}\otimes_RF^0 \to M^{\otimes_RN} \to 0, \,\,\text{and}\,\, \\
0 \to M^{\otimes_RN}\otimes_RF^n \to \cdots \to
M^{\otimes_RN}\otimes_RF^1 \to M^{\otimes_RN}\otimes_RF^0 \to 0. \quad \,\,\,
\end{align*}
Thus, we obtain a flat resolution
\[
0 \to \Ind(F^n) \to \cdots \to \Ind(F^0) \to \Ind(M) \to 0
\]
of $\Ind(M)$, showing that $\fd_{T_R(M)}\Ind(M)<\infty$.
The short exact sequence
\[
0 \to \Ind(M) \to \Ind(R) \to S(R) \to 0
\]
and the fact that $\Ind(R)$ is flat imply $\fd_{T_R(M)} S(R)<\infty$.
\end{prf*}

\begin{lem} \label{PGFTRM in Phi}
If $\fd_RM<\infty$ and $\fd_{\Rop}M<\infty$,
then the inclusion $\PGF(T_R(M)) \subseteq \Phi(\PGF(R))$ holds.
\end{lem}

\begin{prf*}
Let $(X,u) \in \PGF(T_R(M))$.
By Lemma~\ref{projectives in F-Rep}(a), there exists an exact complex of projective $T_R(M)$-modules
\[
Q^\bullet:\
\xymatrix@C=0.5cm{
\cdots \to \Ind(P^{-1}) \ar[r]^-{d^{-1}} &
\Ind(P^{0}) \ar[r]^-{d^0} &
\Ind(P^{1}) \to \cdots
}
\]
with $(X,u)\cong \ker(d^0)$ and each $P^i\in\Proj(R)$, such that
$[E,v]\otimes_{T_R(M)} Q^\bullet$ is exact for every injective $T_R(M)^{\op}$-module $[E,v]$.
Applying the forgetful functor $U$, we obtain an exact complex in $\Mod(R)$:
\[
P^\bullet:\
\xymatrix@C=0.5cm{
\cdots \to
\bigoplus_{i=0}^N M^{\otimes_R i}\otimes_R P^{-1}
\to \bigoplus_{i=0}^N M^{\otimes_R i}\otimes_R P^{0}
\to \cdots.
}
\]

We claim that $M\otimes_R P^\bullet$ is exact.
Indeed, take a flat resolution of $M$ as a right $R$-module,
$0\to F^n\to\cdots\to F^0\to M\to 0.$
Because $M$ satisfies $(\sfT)$,
$\Tor_{\geqslant 1}^R(M,\, \bigoplus_{i=0}^N M^{\otimes_R i}\otimes_R P^j)=0,$
and we obtain an exact sequence of complexes
\[
0 \to F^n\otimes_R P^\bullet \to \cdots \to F^0\otimes_R P^\bullet
\to M\otimes_R P^\bullet \to 0.
\]
Since each $F^i\otimes_R P^\bullet$ is exact, the desired exactness follows.
Hence the diagram
$$
\xymatrix@C=0.7cm@R=0.7cm{
  0 \ar[r]^{}
  & \bigoplus_{i=1}^N(M^{\otimes_Ri}\otimes_RP^{-1}) \ar[d]_{M\otimes d^{-1}}
  \ar[r]^{c_{P^{-1}}}
  & \bigoplus_{i=0}^N(M^{\otimes_Ri}\otimes_RP^{-1})
    \ar[d]_{d^{-1}} \ar[r]^{\qquad \qquad \pi_{P^{-1}}}
  & P^{-1} \ar[d]_{} \ar[r]^{} & 0  \\
  0 \ar[r]^{}
  & \bigoplus_{i=1}^N(M^{\otimes_Ri}\otimes_RP^0) \ar[d]_{M\otimes d^0}
    \ar[r]^{c_{P^0}}
  & \bigoplus_{i=0}^N(M^{\otimes_Ri}\otimes_RP^0)
    \ar[d]_{d^0} \ar[r]^{\qquad \qquad \pi_{P^0}}
  & P^0 \ar[d]_{} \ar[r]^{} & 0  \\
  0 \ar[r]^{}
  & \bigoplus_{i=1}^N(M^{\otimes_Ri}\otimes_RP^1)
  \ar[r]^{c_{P^1}}
  & \bigoplus_{i=0}^N(M^{\otimes_Ri}\otimes_RP^1)
  \ar[r]^{\qquad \qquad \pi_{P^1}}
  & P^1 \ar[r]^{} & 0. }$$
has exact rows and columns, giving a short exact sequence
$0\to M\otimes_R X \xrightarrow{u} X \to \cok(u)\to 0.$

It remains to show that
$I\otimes_RP^\bullet$ is exact for each injective $\Rop$-module $I$,
where $P^\bullet$ is the third non-zero column in the above commutative diagram.
If we are done, then $\cok(u)\in\PGF(R)$.
Indeed, there exists a split monomorphism $I \to \prod R^+$,
which induces the split monomorphism of $T_R(M)^{\op}$-modules
$$\xymatrix{0 \to S(I) \to S(\prod R^+) \cong \prod S(R)^+}.$$
Since $\fd_{T_R(M)}S(R)<\infty$ by Lemma \ref{dim of S(R)},
one can easily see that $\id_{T_R(M)^{\op}}S(R)^+<\infty$,
and hence $\id_{T_R(M)^{\op}}S(I)<\infty$.
Note that $$S(I)\otimes_{T_R(M)^{\op}}Q^\bullet \cong I\otimes_RC(Q^\bullet)$$
by Corollary \ref{cor fof tensor products}.
Since the left-hand complex is exact and $C(Q^\bullet)=P^\bullet$
(by Remark \ref{adjoint pairs}),
we obtain that $I\otimes_R P^\bullet$ is exact.
\end{prf*}

\begin{lem} \label{Phi in PGFTRM}
If $\pd_RM<\infty$ and $\fd_{\Rop}M<\infty$,
then the inclusion $\Phi(\PGF(R)) \subseteq \PGF(T_R(M))$ holds.
\end{lem}
\begin{prf*}
Let $(X,u) \in \Phi(\PGF(R))$, and consider the short exact sequence
\[
\delta:\ 0 \to M\otimes_RX \overset{u \,}\longrightarrow X
\overset{\pi \,}\longrightarrow \cok(u) \to 0
\]
in $\Mod(R)$.
Since $\cok(u) \in \PGF(R)$, there exists an exact sequence
$$\xymatrix@C=0.5cm{
P^\bullet: \cdots \to P^{-1} \ar[r]^{\qquad \quad d^{-1} \,} & P^0
\ar[r]^{ d^0 \quad \,\,} & P^1 \to \cdots}$$
of projective $R$-modules with $\cok(u) \cong \kernel(P^0 \to P^1)$,
such that it remains exact after applying the functor
$E\otimes_R-$ for every $\Rop$-module $E$ with finite injective dimension.
Note that $M$ is admissible by \prpcite[2.7]{2025DLST}.
So, the usual construction (cf. \prpcite[2.11]{2025DLST}) yields an exact complex of projective $T_R(M)$-modules
$$\Delta:\ \xymatrix@C=0.5cm{
\cdots \to \Ind(P^{-1}) \ar[r] & \Ind(P^{0})
\ar[r] & \Ind(P^{1}) \to \cdots}$$ with $(X,u)\is\ker(d^0)$.

We must show that $\Coind(E)\otimes_{T_R(M)} \Delta$ is exact for all injective $\Rop$-modules $E$.
For each $1 \leqslant i \leqslant N$,
set $(M^{\otimes_R i},E):=\Hom_{\Rop}(M^{\otimes_R i},E)$.
There are short exact sequences
\[
0\to S((M^{\otimes_R i},E)) \to \Coind((M^{\otimes_R i},E))
\to \Coind((M^{\otimes_R(i+1)},E)) \to 0,
\]
and when $i=N$ we have
$S((M^{\otimes_R N},E)) \cong \Coind((M^{\otimes_R N},E)).$
Thus, it suffices to show that
\[
S((M^{\otimes_R i},E))\otimes_{T_R(M)} \Delta
\]
is exact for all $i$.
If we are done, then the isomorphism
$$S((M^{\otimes_RN},E))\otimes_{T_R(M)}\Delta
\cong \Coind((M^{\otimes_RN},E))\otimes_{T_R(M)}\Delta$$
implies that the complex
$\Coind((M^{\otimes_RN-1},E))\otimes_{T_R(M)}\Delta$ is exact.
By induction, one can easily see that
the complex $\Coind(E)\otimes_{T_R(M)}\Delta$ is also exact.
Indeed, by Lemma \ref{lemma4.5 for (T)}, we have $\fd_RM^{\otimes_Ri}<\infty$.
This implies that $\id_{\Rop}(M^{\otimes_Ri},E)<\infty$
by \lemcite[2.2]{GFGFD}. However,
$$S((M^{\otimes_Ri},E))\otimes_{T_R(M)}\Delta
\cong (M^{\otimes_Ri},E)\otimes_RP^\bullet$$
by Corollary \ref{cor fof tensor products}.
Thus, the complexes $S(\Hom_{\Rop}(M^{\otimes_Ri},E))\otimes_{T_R(M)}\Delta$
are exact, as desired.
\end{prf*}

Combining Lemmas~\ref{PGFTRM in Phi} and \ref{Phi in PGFTRM} yields Theorem~\ref{THM PGF com}.
Let $X$ be an $R$-module.
Note that $X=C \circ \Ind(X)$ by \ref{adjoint pairs}.
This implies that $X \in \PGF(R)$ if and only if $\Ind(X) \in \Phi(\PGF(R))$.
Hence, the following result is an immediate consequence of Theorem~\ref{THM PGF com}.

\begin{cor}\label{PGF IND}
If $\pd_RM<\infty$ and $\fd_{\Rop}M<\infty$,
then $X$ is a PGF $R$-module if and only if $\Ind(X)$ is a PGF $T_R(M)$-module.
\end{cor}

\section{Proof of Theorem \ref{THM GF com}}
\noindent
This section is devoted to the proof of Theorem~\ref{THM GF com}. Throughout, we assume that the $R$-bimodule $M$ satisfies condition~$(\sfT)$. We begin by recalling the notion of Gorenstein flat modules.

\begin{dfn} \label{def of GF} \cite{GF1993}
An $R$-module $X$ is called \emph{Gorenstein flat} if there is an exact complex
$$\xymatrix@C=0.5cm{
F^\bullet: \cdots \to F^{-1} \ar[r]^{\qquad \quad d^{-1} \,} & F^0
\ar[r]^{\,\, d^0 \qquad } & F^1 \to \cdots,}$$
of flat $R$-modules such that $X \cong \ker(F^0 \to F^1)$ and the complex
$E\otimes_R F^\bullet$ is exact for every injective $\Rop$-module $E$.
We denote by $\GF(R)$ the class of all Gorenstein flat $R$-modules.
\end{dfn}

\begin{rmk}\label{id finite}
By \lemcite[2.1]{GFGFD}, an $R$-module $X$ is Gorenstein flat if and only if there exists a flat resolution $F^\bullet$ as above such that $E\otimes_R F^\bullet$ remains exact for every $\Rop$-module $E$ with finite injective dimension.
\end{rmk}

\begin{lem} \label{X in GFTRM}
Suppose that $\fd_RM<\infty$ and $\fd_{\Rop}M<\infty$.
If $X$ is a Gorenstein flat $R$-module,
then $\Ind(X)$ is a Gorenstein flat $T_R(M)$-module.
\end{lem}
\begin{prf*}
Suppose that $X \in \GF(R)$. Then there exists an exact complex
$$\xymatrix@C=0.5cm{
F^\bullet: \cdots \to F^{-1} \ar[r]^{\qquad \quad d^{-1} \,} & F^0
\ar[r]^{\,\,\, d^0 \qquad } & F^1 \to \cdots,}$$
of flat $R$-modules with $X \cong \kernel(F^0 \to F^1)$,
such that it remains exact after applying the functor
$E\otimes_R-$ for every $\Rop$-module $E$ with finite injective dimension.
For each $1 \leqslant i \leqslant N$,
since $\fd_{\Rop}M^{\otimes_Ri}<\infty$ by Lemma \ref{lemma4.5 for (T)},
$M^{\otimes_Ri}\otimes_RF^\bullet$ is exact by \lemcite[2.3]{GCOND}.
So we get the exact complex of flat $T_R(M)$-modules
$$\Delta:\ \xymatrix@C=0.5cm{
\cdots \to \Ind(F^{-1}) \ar[r] & \Ind(F^{0})
\ar[r] & \Ind(F^{1}) \to \cdots}$$ with $\Ind(X) \is\ker(d^0)$.
To complete the proof, it remains to show that the complex
$\Coind(E)\otimes_{T_R(M)}\Delta$ is exact for each $E \in \Inj(\Rop)$.
Recall the second paragraph of the proof of Lemma \ref{Phi in PGFTRM}.
For the same reasons, it suffices to show that
$S(\Hom_{\Rop}(M^{\otimes_Ri},E))\otimes_{T_R(M)}\Delta$ is exact,
which is guaranteed by $\fd_RM<\infty$.
\end{prf*}

\begin{prp} \label{GFTRM in X}
Suppose that $\fd_RM<\infty$ and $\fd_{\Rop}M<\infty$.
If $(X,u)$ is a Gorenstein flat $T_R(M)$-module,
then $\cok(u)$ is a Gorenstein flat $R$-module.
\end{prp}
\begin{prf*}
Let $(X,u)$ be in $\GF(T_R(M))$.
Then, by Lemma \ref{projectives in F-Rep}(c),
there exists an exact complex of flat $T_R(M)$-modules
$$\Delta:\ \xymatrix@C=0.5cm{
\cdots \to (G^{-1}, h^{-1}) \ar[r]^-{d^{-1}} & (G^{0}, h^{0})
\ar[r]^-{d^0} & (G^{1}, h^{1}) \to \cdots}$$
with $(X,u)\is\ker(d^0)$ and each $\cok(h^{i}) \in \Flat(R)$,
such that $[E,v]\otimes_{T_R(M)}\Delta$ is exact
for every $T_R(M)^{\op}$-module $[E,v]$ with finite injective dimension.
Since $\fd_{{T_R(M)}^{\op}}S(R)<\infty$ by Lemma \ref{dim of S(R)},
the sequence $S(R)\otimes_{T_R(M)}\Delta$ is exact by \lemcite[2.3]{GCOND}.
Note that the isomorphisms
\begin{align*}
&S(R)\otimes_{T_R(M)}(X,u) \cong R\otimes_R\cok(u) \cong \cok(u) \text{\,\, and} \\
&S(R)\otimes_{T_R(M)}(G^i, h^i) \cong R\otimes_R\cok(h^i) \cong \cok(h^i)
\end{align*}
hold by Corollary \ref{cor fof tensor products}.
Hence, there exists an exact complex of flat $R$-modules
$$G^\bullet: \cdots \to \cok(h^{-1}) \to  \cok(h^{0}) \to  \cok(h^{1}) \to \cdots$$
with $\cok(u)\is \ker(\cok(h^{0}) \to \cok(h^{1}))$.
Let $I$ be an injective $\Rop$-module.
Recall the third paragraph of the proof of Lemma \ref{PGFTRM in Phi}.
For the same reasons, one can easily see that
the sequence $I\otimes_RG^\bullet$ is exact, since $\fd_RM<\infty$ by assumption.
This completes the proof.
\end{prf*}

The following result is an immediate consequence of Proposition \ref{GFTRM in X}.
Now, Theorem \ref{THM GF com}
follows from Lemma \ref{X in GFTRM} and Corollary \ref{GFTRM in IndX}.

\begin{cor} \label{GFTRM in IndX}
Suppose that $\fd_RM<\infty$ and $\fd_{\Rop}M<\infty$.
If $\Ind(X)$ is a Gorenstein flat $T_R(M)$-module,
then $X$ is a Gorenstein flat $R$-module.
\end{cor}

\begin{rmk} \label{rmk for GFTRM}
It is instructive to compare our Theorem \ref{THM GF com} with a result derived from \thmcite[3.8]{2025DLST}.
The latter states that if
$R$ is right coherent, and if $M$ is flat as an $R$-module
and finitely presented with finite projective dimension as an $\Rop$-module,
then an $R$-module $X$ is Gorenstein flat if and only if
$\Ind(X)$ is a Gorenstein flat $T_R(M)$-module.
Theorem \ref{THM GF com} achieves a similar characterization under a set of
significantly weaker hypotheses.
\end{rmk}

The following example is due to Cui, Rong and Zhang \exacite[4.7]{2025ZhangPu},
which shows that there exists a bimodule $M$ over some algebra satisfying
all the conditions in Theorems \ref{THM PGF com} and \ref{THM GF com}.

\begin{exa} \label{example of qviver}
Let $k$ be a field and $kQ$ the path algebra
associated to the quiver
\begin{align*}
\xymatrix{
Q:\  1 \ar[r]_{} & 2 \ar[r]_{} & 3 \ar[r]_{}
  & \cdots \ar[r]_{} & n \ar@/_3ex/[llll]_{}.}
\end{align*}
Suppose that $J$ is the ideal of $kQ$ generated by all the arrows.
Then $R=kQ/J^h$ is a self-injective algebra for $2\leqslant h \leqslant n$.
Denote by $e_i$ the idempotent element corresponding to the vertex $i$.
Then one has $e_jRe_i = 0$ whenever
$1\leqslant i< j\leqslant n$ and $j-i \geqslant h$.
Let $M = Re_i\otimes_k e_jR$. Then $M$ is an $R$-bimodule and projective on both sides, and $M\otimes_RM\cong Re_i\otimes_k (e_jR\otimes_RRe_i)\otimes_k e_jR=0$.
\end{exa}

\section{Applications}
\label{Applications}
\noindent
In this section, we give some applications of our results to trivial ring extensions and Morita context rings.

\subsection{The trivial extension of rings}
\label{The trivial extension of rings}
\noindent
Let $R$ be a ring and $M$ an $R$-bimodule.
There exists a ring $R\ltimes M$,
where the addition is componentwise
and the multiplication is given by
$(r_1, m_1)(r_2, m_2) = (r_1r_2,r_1m_2 + m_1r_2)$ for $r_1, r_2 \in R$ and $m_1, m_2 \in M$.
This ring is called the \emph{trivial extension} of
the ring $R$ by the $R$-bimodule $M$; see \cite{TRIEXT1975} and \cite{TRIEXT1971}.

Suppose that the $R$-bimodule $M$ is $1$-nilpotent,
that is, $M\otimes_RM=0$.
Then it is straightforward to verify that the tensor ring
$T_R(M)$ coincides with the trivial ring extension $R\ltimes M$.
One can immediately get the following results by
Theorems \ref{THM PGF com} and \ref{THM GF com}, respectively.

\begin{cor} \label{PGF in trivial extension}
Suppose that $M$ is a $1$-nilpotent $R$-bimodule.
If $\Tor_{\geqslant 1}^R(M,M\otimes_RP)=0$ for each $P\in \Proj(R)$,
$\pd_{R}M<\infty$ and $\fd_{\Rop}M<\infty$,
then there is an equality $\PGF(R\ltimes M) = \Phi(\PGF(R))$.
\end{cor}

\begin{cor} \label{GF in trivial extension}
Suppose that $M$ is a $1$-nilpotent $R$-bimodule.
If $\Tor_{\geqslant 1}^R(M,M\otimes_RP)=0$ for each $P\in \Proj(R)$,
$\fd_RM<\infty$ and $\fd_{\Rop}M<\infty$,
then an $R$-module $X$ is Gorenstein flat if and only if
$\Ind(X)$ is a Gorenstein flat $R\ltimes M$-module.
\end{cor}

\begin{rmk} \label{GP in trivial extension}
Recall Lemma \ref{dim of S(R)}. Under the assumption that
$M$ is a $1$-nilpotent $R$-bimodule such that
$\Tor_{\geqslant 1}^R(M,M\otimes_RP)=0$ for each $P\in \Proj(R)$,
we may replace the assumptions in \thmcite[2.11]{2025PGF}, namely
$\fd_{T_R(M)}S(R)<\infty$ and $\fd_{{T_R(M)}^{\op}}S(R)<\infty$,
by the conditions $\fd_{R}M<\infty$ and $\fd_{\Rop}M<\infty$.
Similarly for the assumptions in \corcite[4.8(2)]{TriExtGPMao}.
\end{rmk}

\subsection{Morita context rings}
\label{Morita context rings}
Let $A$ and $B$ be two rings,
and let $_AV_B$ and $_BU_A$ be two bimodules,
$\phi : U\otimes_AV \to B$ a homomorphism of $B$-bimodules,
and $\psi : V\otimes_BU \to A$ a homomorphism of $A$-bimodules.
We assume further that $\phi(u\otimes v)u'=u\psi(v \otimes u')$
and $v\phi(u\otimes v')=\psi(v \otimes u)v'$
for all $u, u' \in U$ and $v, v' \in V$.
Associated with a \emph{Morita context} $(A,B,U,V,\phi,\psi)$,
there exists a \emph{Morita context ring}
\[\Lambda_{(\phi,\psi)}=\begin{pmatrix}A & V \\ U & B\end{pmatrix},\]
where the addition of elements is componentwise and multiplication is given by
\[\begin{pmatrix}a & v \\ u & b\end{pmatrix}\cdot
\begin{pmatrix}a' & v' \\ u' & b'\end{pmatrix}=
\begin{pmatrix}aa'+\psi(v \otimes u') & av'+vb' \\
ua'+bu' & bb'+ \phi(u \otimes v') \end{pmatrix};\]
see \cite{BASS1968, MORITA1958} for more details.
Following \thmcite[1.5]{GF1982}, one can view a $\Lambda_{(\phi,\psi)}$-module
as a quadruple $(X,Y,f,g)$ with $X \in \Mod(A)$, $Y \in \Mod(B)$,
$f \in \Hom_B(U\otimes_AX,Y)$, and $g \in \Hom_A(V\otimes_BY,X)$. Also, one can view a $\Lambda_{(\phi,\psi)}^{\op}$-module as a quadruple
$[W,N,s,t]$ with $W \in \Mod(A^{\op})$, $N \in \Mod(B^{\op})$,
$s \in \Hom_{B^{\op}}(N,\Hom_{A^{\op}}(U,W))$, and $t \in \Hom_{A^{\op}}(W,\Hom_{B^{\op}}(V,N))$.

It follows from \prpcite[2.5]{GFARTIN} that Morita context rings are trivial ring extensions
whenever both $\phi$ and $\psi$ are zero.
More precisely, consider the Morita context ring $\Lambda_{(0,0)}$.
There exists an isomorphism of rings:
$$\Lambda_{(0,0)}\overset{\cong}\longrightarrow (A\times B)\ltimes (U \oplus V)
\,\,\text{via}\,\, \begin{pmatrix}a & v \\ u & b\end{pmatrix} \mapsto
((a,b),(u,v)).$$
Thus, there exists an isomorphic functor
$$\mu: \Mod(\Lambda_{(0,0)}) \to \Mod((A\times B)\ltimes (U \oplus V))\ \mathrm{via}\
(X,Y,f,g) \mapsto ((X,Y),(g,f)),$$
where $(g,f)$ is from
$(U \oplus V)\otimes_{A\times B}(X, Y) \cong (V\otimes_BY,U\otimes_AX)$ to $(X, Y)$.

Note that $(U \oplus V)\otimes_{A\times B}(U \oplus V) \cong (U\otimes_AV) \oplus (V\otimes_BU)$.
Then the $A\times B$-bimodule $U \oplus V$ is $1$-nilpotent if and only if $U\otimes_AV=0=V\otimes_BU$.
Thus, we obtain the following result by Corollary \ref{PGF in trivial extension}.

\begin{cor} \label{PGP in morita ring}
Let $\Lambda_{(0,0)}$ be a Morita context ring with $U\otimes_AV=0=V\otimes_BU$.
Suppose that $\Tor^B_{\geqslant 1}(V, U\otimes_AP_1)=0=\Tor^A_{\geqslant 1}(U, V\otimes_BP_2)$ for each $P_1 \in \Proj(A)$ and $P_2 \in \Proj(B)$.
If $\pd_{B}U<\infty$, $\fd_{A^{\op}}U<\infty$,
$\pd_{A}V<\infty$ and $\fd_{B^{\op}}V< \infty$,
then a $\Lambda_{(0,0)}$-module $(X,Y,f,g)\in\PGF(\Lambda_{(0,0)})$
if and only if both $f$ and $g$ are monomorphisms,
$\cok(f)\in\PGF(B)$ and $\cok(g)\in\PGF(A)$.
\end{cor}
\begin{prf*}
By the assumption, one can check that
$\Tor^{A\times B}_{\geqslant 1}
(U \oplus V, (U \oplus V)\otimes_{A\times B}(P_1,P_2))=0$ for each
$(P_1,P_2) \in \Proj(A\times B)$, and both $\pd_{A\times B}(U \oplus V)$ and $\fd_{{(A\times B)}^{\op}}(U \oplus V)$ are finite.
Thus, it follows from Corollary \ref{PGF in trivial extension}
that there is an equality
$\PGF(\Lambda_{(0,0)})
\cong \PGF((A\times B)\ltimes (U \oplus V))
=\Phi(\PGF(A\times B)),$
where an object $((X,Y),(g,f))\in\Mod((A\times B)\ltimes (U \oplus V))$ is in
$\Phi(\PGF(A\times B))$ if and only if $(g,f)$
is a monomorphism and $\cok(g,f)$ is in $\PGF(A\times B)$.
This yields that $(X,Y,f,g)\in\PGF(\Lambda_{(0,0)})$ if and only if
both $f$ and $g$ are monomorphisms, and $\cok(f)\in\PGF(B)$ and $\cok(g)\in\PGF(A)$.
\end{prf*}

\begin{rmk} \label{rmk for PGP in morita ring}
Under the assumption that $U\otimes_AV=0=V\otimes_BU$
and $\Tor^B_{\geqslant 1}(V, U\otimes_AP_1)=0=\Tor^A_{\geqslant 1}(U, V\otimes_BP_2)$
for each $P_1 \in \Proj(A)$ and $P_2 \in \Proj(B)$,
we may replace the assumptions in \prpcite[3.4]{2025PGF}
that $\fd_{\Lambda_{(0,0)}}(A,B,0,0)<\infty$ and
$\fd_{{\Lambda_{(0,0)}}^{\op}}(A,B,0,0)<\infty$
with the conditions that $\pd_{B}U<\infty$, $\fd_{A^{\op}}U<\infty$,
$\pd_{A}V<\infty$ and $\fd_{B^{\op}}V< \infty$.
\end{rmk}

Similarly, by Corollary \ref{GF in trivial extension}, one obtains the following result.
Also, under the additional assumption that $U\otimes_AV=0=V\otimes_BU$
and $\Tor^B_{\geqslant 1}(V, U\otimes_AQ)=0=\Tor^A_{\geqslant 1}(U, V\otimes_BP)$
for each $Q \in \Proj(A)$ and $P \in \Proj(B)$,
we can remove the hypothesis in \thmcite[2.5]{Mao24GF} that
$\fd_{\Lambda_{(0,0)}}(A,0,0,0)<\infty$, $\fd_{{\Lambda_{(0,0)}}^{\op}}(A,0,0,0)<\infty$,
$\fd_{\Lambda_{(0,0)}}(0,B,0,0)<\infty$ and $\fd_{{\Lambda_{(0,0)}}^{\op}}(0,B,0,0)<\infty$.

\begin{cor} \label{GF in morita ring}
Let $\Lambda_{(0,0)}$ be a Morita context ring with $U\otimes_AV=0=V\otimes_BU$.
\begin{prt}
\item
Suppose that $\Tor^A_{\geqslant 1}(U, V\otimes_BP)=0$ for each $P \in \Proj(B)$.
If $\fd_{A}V<\infty$ and $\fd_{A^{\op}}U<\infty$,
then an $A$-module $X$ is Gorenstein flat if and only if
$(X,U\otimes_AX,1,0)$ is a Gorenstein flat $\Lambda_{(0,0)}$-module;
\item
Suppose that $\Tor^B_{\geqslant 1}(V, U\otimes_AQ)=0$ for each $Q \in \Proj(A)$.
If $\fd_{B}U<\infty$ and $\fd_{B^{\op}}V<\infty$,
then a $B$-module $Y$ is Gorenstein flat if and only if
$(V\otimes_BY,Y,0,1)$ is a Gorenstein flat $\Lambda_{(0,0)}$-module.
\end{prt}
\end{cor}

We end this paper by presenting the following results for
the formal lower triangular matrix ring $T=\begin{pmatrix}A & 0 \\ U & B\end{pmatrix}$.

\begin{cor} \label{PGF over T} \corcite[3.7]{2025PGF}
If $\pd_{B}U<\infty$ and $\fd_{A^{\op}}U<\infty$,
then a $T$-module $(X,Y,f,0)\in\PGF(T)$
if and only if $f$ is a monomorphism,
$X \in\PGF(A)$ and $\cok(f)\in\PGF(B)$.
\end{cor}

\begin{cor} \label{GF over T}
If $\fd_{B}U<\infty$ and $\fd_{A^{\op}}U<\infty$,
then an $A$-module $X$ is Gorenstein flat if and only if
$(X,U\otimes_AX,1,0)$ is a Gorenstein flat $T$-module.
\end{cor}


\bibliographystyle{amsplain-nodash}

\def\cprime{$'$}
  \providecommand{\arxiv}[2][AC]{\mbox{\href{http://arxiv.org/abs/#2}{\sf
  arXiv:#2 [math.#1]}}}
  \providecommand{\oldarxiv}[2][AC]{\mbox{\href{http://arxiv.org/abs/math/#2}{\sf
  arXiv:math/#2
  [math.#1]}}}\providecommand{\MR}[1]{\mbox{\href{http://www.ams.org/mathscinet-getitem?mr=#1}{#1}}}
  \renewcommand{\MR}[1]{\mbox{\href{http://www.ams.org/mathscinet-getitem?mr=#1}{#1}}}
\providecommand{\bysame}{\leavevmode\hbox to3em{\hrulefill}\thinspace}
\providecommand{\MR}{\relax\ifhmode\unskip\space\fi MR }
\providecommand{\MRhref}[2]{%
  \href{http://www.ams.org/mathscinet-getitem?mr=#1}{#2}
}
\providecommand{\href}[2]{#2}

\end{document}